\documentclass[12pt]{article}      
\usepackage{amsmath,amssymb,amsthm, amscd, graphicx, verbatim, setspace} 
\usepackage{setspace}
\usepackage{mathtools}
\usepackage{tabstackengine}

\theoremstyle{plain}
\newtheorem{thm}{Theorem}

\usepackage{algorithm}
\usepackage{algpseudocode}
\usepackage{amsthm}
\usepackage{amsmath}

\title{Improved lower bound on generalized Erdos-Ginzburg-Ziv constants}
\date{}
\author{Jesse Geneson (PSU)}
\begin{document}
\maketitle

\begin{abstract}
If $G$ is a finite Abelian group, define $s_{k}(G)$ to be the minimal $m$ such that a sequence of $m$ elements in $G$ always contains a $k$-element subsequence which sums to zero. Recently Bitz et al. proved that if $n = exp(G)$, then $s_{2n}(C_{n}^{r}) > \frac{n}{2}[\frac{5}{4}-O(n^{-\frac{3}{2}})]^{r}$ and $s_{k n}(C_{n}^{r}) > \frac{k n}{4} [1+\frac{1}{e k}-O(\frac{1}{n})]^{r}$ for $k > 2$. In this note, we sharpen their general bound by showing that $s_{k n}(C_{n}^{r}) > \frac{k n}{4} [1+\frac{(k-1)^{(k-1)}}{k^k}-O(\frac{1}{n})]^{r}$ for $k > 2$.
\end{abstract}

\section{Lower bound}

The function $s_{k n}(G)$ is known as the $k^{th}$ generalized Erdos-Ginzburg-Ziv constant of $G$. The first result about these constants was proved in \cite{egz}, and there has been a long history of results since then \cite{cl, eg, ga, ha, ku, re, ro}, which were detailed in \cite{bgh}. Our improvement is the theorem below.

\begin{thm}\label{main}
For $k > 2$, we have $s_{k n}(C_{n}^{r}) > \frac{k n}{4} [1+\frac{(k-1)^{(k-1)}}{k^k}-O(\frac{1}{n})]^{r}$.
\end{thm}

To prove this result, we use the bound of Sondow et al. \cite{so, sz} for binomial coefficients of the form $\binom{k n}{n}$. Specifically they proved the bounds $\frac{1}{4(k-1)n}[\frac{k^k}{(k-1)^{(k-1)}}]^{n} < \binom{k n}{n} < [\frac{k^k}{(k-1)^{(k-1)}}]^{n}$ for $n$ a positive integer and $k \geq 2$ a real number. The proof of Theorem \ref{main} is nearly identical to the proof in \cite{bgh}, but just with this sharper bound on the binomial coefficients. 

As in \cite{bgh}, define $N = \frac{k n A^{r}}{4}$ for $A$ to be chosen at the end of the proof. Randomly pick a sequence $X$ of $N$ vectors in $\left\{0,1\right\}^{r}$ by letting each coordinate be $1$ with probability $q$, and let $Z$ be the number of subsequences of length $k n$ in $X$ that sum to $0$. We show that $E[Z] < 1$ with $A = 1+\frac{(k-1)^{(k-1)}}{k^k}-O(\frac{1}{n})$.

First we calculate the probability $Q$ that a given coordinate sums to $0$, which is the sum of the probabilities $P_{i n}$ that the coordinate sums to $i n$ for $0 \leq i \leq k$, which are equal to $P_{i n} = \binom{k n}{i n} q^{i n} (1-q)^{(k-i)n}$. We want the terms $P_{0}$ and $P_{n}$ to dominate, so we set $(1-q)^{k n} = \binom{k n}{n}(1-q)^{(k-1)n} q^{n}$. Combining the bounds of Sondow et al. with this equality, we obtain $\frac{1}{4(k-1)n}[\frac{k^k}{(k-1)^{(k-1)}}]^{n} (1-q)^{(k-1)n} q^{n} < (1-q)^{k n} < [\frac{k^k}{(k-1)^{(k-1)}}]^{n} (1-q)^{(k-1)n} q^{n}$. This implies that $(1-o(1))[\frac{k^k}{(k-1)^{(k-1)}}] q < 1-q < [\frac{k^k}{(k-1)^{(k-1)}}] q$, and thus that $\frac{1}{1+\frac{k^k}{(k-1)^{(k-1)}}} < q < \frac{1}{1+(1-o(1))\frac{k^k}{(k-1)^{(k-1)}}}$.

$P_{0}$ and $P_{n}$ are approximately equal in this range, so that $Q < (k+1)(1-q)^{k n}$ and $E[Z] = \binom{N}{k n} Q^{r} < (\frac{4N}{k n})^{k n} (k+1)^{r} (1-q)^{k n r}$. We want $E[Z] < 1$, so $A < \frac{1}{(k+1)^{1/kn}(1-q)}$ will suffice, and thus we may let $A = 1+\frac{(k-1)^{(k-1)}}{k^k}-O(\frac{1}{n})$.

\end{document}